\documentstyle{article}

\begin{document}
\newcommand{\nt}{\noindent}

\begin{center}

{\bf Using Symmetry to Count Rational Curves} 

\medskip

Aaron Bertram\footnote{Supported in part by NSF Research Grant
DMS-9970412}

\medskip

{\it Dedicated to Herb Clemens} 

\medskip

\end{center}

The recent ``close encounter'' between enumerative algebraic
geometry and theoretical physics has resulted in many new
applications and new techniques for
counting algebraic curves on a complex
projective manifold. For example, string theorists
demonstrated that generating functions built from counting
curves can have completely unexpected
relationships with other geometric
constructions, via mirror symmetry \cite{Can}. In this paper, I want
to focus on generating functions built from counting rational curves,
and another insight inspired by physicists -- the existence
of hidden symmetries in the generating functions themselves.

\medskip

By ``counting rational curves,'' I mean producing numbers such as:

\medskip

$\bullet$ the number of lines in ${\bf P}^n$ through two general points
$(1)$

\medskip

$\bullet$ the number of lines in ${\bf P}^3$ meeting $4$ general lines
$(2)$

\medskip

$\bullet$ the number of rational quartics
in ${\bf P}^2$ through $11$ general points $(620)$

\medskip

$\bullet$ the expected number of twisted cubic curves on a general quintic
in ${\bf P}^6$ passing through $2$ general points $(2,088,000)$

\medskip

The first number is basic, the second is ``classical'',
but the last two were not known until quite recently. There are now
several methods for computing them \cite{Ell,Gath,Ran}, but the goal
of this paper is to show how easily such computations follow from 
considerations of symmetry. More precisely, I will
review WDVV symmetry and the reconstruction theorem of Kontsevich and Manin,
then introduce an analogous family of new symmetries and 
corresponding reconstruction
theorem recently discovered in joint work together with Holger Kley
\cite{BerK}. 

\medskip

Let $X$ be a complex projective manifold, and define:

\medskip

$V = H^*(X,{\bf Q})$ with 
the Poincar\'e pairing $\gamma_1\otimes \gamma_2 \mapsto \int_X \gamma_1
\wedge \gamma_2$

\medskip

$C_{eff} \subset H_2(X,{\bf Z})$, the set of positive linear
combinations of classes $f_*[C]$, where
$C$ is a Riemann surface and $f:C \rightarrow X$ is a morphism.
Since $C_{eff}$ has the property that for 
each $\beta \in C_{eff}$:
$$C_{eff} \cap (\beta - C_{eff}) = \{\alpha \in C_{eff} | \beta - \alpha \in C_{eff}\}
\
\mbox{is finite}$$
it follows that the ``Novikov ring'' ${\bf Q}[[q]]$, which is the set
of  (infinite) sums
$\sum_{\beta\in C_{eff}} a_\beta q^\beta$ with $a_\beta\in {\bf Q}$,
has a well-defined multiplication with $q^\alpha q^\beta = q^{\alpha +
\beta}$.

\medskip

All the generating functions will be defined in terms of $V$ and
${\bf Q}[[q]]$.

\medskip

\nt {\bf 1. WDVV Symmetry}. The Gromov-Witten invariant
 $\langle
\gamma_1,...,\gamma_n\rangle_\beta$, for given 
$\beta\in C_{eff}$ and 
$\gamma_1,...,\gamma_n \in V$, counts the expected number of 
rational curves of class $\beta$ ``meeting'' the classes $\gamma_i$
(as in the earlier examples). We will show below how this definition
is made rigorous by means of an explicit 
symmetric element $c^\beta_{0,n}(X)\in
V^{\otimes n} \cong H^*(X^n,{\bf Q})$ with the definition:
$$\langle \gamma_1,...,\gamma_n\rangle_\beta := \int_{X^n} c^\beta_{0,n}(X)
\wedge \gamma_1\otimes ...\otimes \gamma_n$$ 
\nt Think of $Sym^*V$
as the permutation-invariant
sub-algebra of the tensor algebra, and define:
$$c_0(X) := \sum_{\beta}\sum_{n \ge 3} c^\beta_{0,n}(X)q^\beta \in 
Sym^*V[[q]] = Sym^*V\otimes_{\bf Q}{\bf Q}[[q]]$$
This starts with $n=3$ because of the impossibility of finding a good
definition for the $c^0_{0,n}(X)$ for $n < 3$ (but all
$c^\beta_{0,n}(X)$ do have good definitions when $\beta \ne 0$).

\medskip

So far, $c_0(X)$ is just a convenient means of packaging the Gromov-Witten
invariants, but a key insight coming from string theory is the
hidden:

\medskip

\nt {\bf WDVV Symmetry:} For the contraction ``$\rfloor$'' defined below,
$$c_0(X)\rfloor c_0(X) \in Sym^*V\otimes V^{\otimes 4}[[q]]$$
is invariant under the permutation action of the symmetric group $S_4$ on
$V^{\otimes 4}$.

\medskip

The contraction is defined by partially
desymmetrizing, thinking of:
$$c_0(X) \in Sym^*V\otimes V^{\otimes 3}[[q]]$$
then isolating the last copy of $V$ from each $V^{\otimes 3}$ and
contracting them via:
$$V[[q]]\otimes V[[q]] \rightarrow {\bf Q}[[q]];\ \ \gamma_1q^{\alpha_1}
\otimes \gamma_2q^{\alpha_2} \mapsto \left(\int_X \gamma_1\wedge
\gamma_2\right) q^{\alpha_1+\alpha_2},$$
and finally partially symmetrizing via the map:
$Sym^*V \otimes Sym^*V \rightarrow Sym^*V$.

\medskip

\nt {\bf Example:} The coefficient of $q^0$ in $c_0(X)$ is the diagonal
$\Delta_3 \in V^{\otimes 3}$, and:
$$\Delta_3 \rfloor \Delta_3 = \Delta_4 \in V^{\otimes 4}$$
is symmetric. On the other hand, an arbitrary symmetric element of
$V^{\otimes 3}$ does not result in a symmetric element of $V^{\otimes 4}$.
For instance, if $X = {\bf P}^1$, then:
$$(1^{\otimes 3} + H^{\otimes 3})\rfloor (1^{\otimes 3} + H^{\otimes 3})
= 1^{\otimes 2}\otimes H^{\otimes 2} + H^{\otimes 2}\otimes 1^{\otimes 2}$$
is symmetric under the subgroup of $S_4$ generated by $(1 2),(3 4)$ and $(1
3)(2 4)$ 
(These symmetries always hold.) The extra $(2\ 3)$ symmetry is the
novel one.

\medskip

At this point, a word is in order on the construction of the $c^\beta_{0,n}(X)$.

\medskip
 
\nt {\bf Definition:} A pointed rational curve $p_1,...,p_n \in C$ is {\bf stable}
if $C$ is a tree of smooth rational curves with nodes, the $p_1,...,p_n\in C$
are distinct smooth points, and the automorphism group Aut$(C;p_i)$ of $C$, fixing
the $p_i$, is finite.

\medskip

\nt {\bf Theorem:} (Mumford-Knudsen) If $n \ge 3$, a fine moduli space of stable
pointed rational curves exists. It is projective, smooth of dimension $n-3$,
denoted
$\overline M_{0,n}$.

\medskip

\nt {\bf Examples:} 
$\overline M_{0,3} =$ point,
$\overline M_{0,4}
\cong {\bf P}^1$.

\medskip

Inspired by this, Kontsevich-Manin made the following definition
 \cite{KonM}:

\medskip

\nt {\bf Definition:} A map $f:(C;p_1,...,p_n) \rightarrow X$ from a tree of smooth rational curves
with distinct smooth points is {\bf stable} if the automorphism group
Aut$(C;p_i;f)$ of  $C$, fixing the $p_i$ and the map $f$, is
finite.

\medskip

\nt {\bf Examples:} All maps with finite fibers are stable. A constant map is
stable if and only if the pointed curve is stable.

\medskip

\nt {\bf Theorem:} (Kontsevich-Manin) A proper, Deligne-Mumfod stack
of stable maps, denoted $\overline M_{0,n}(X,\beta)$, exists for each
fixed $\beta$. It is smooth of dimension
$\mbox{dim}(X) - K_X \cap \beta + n - 3$ when $X$ is nice enough
(e.g. a homogeneous space).

\medskip

Moreover, the moduli space comes equipped with ``structure'' maps:
$$\begin{array}{ccc}\overline M_{0,n}(X,\beta) & \stackrel {ev}\rightarrow & X^n \\
\pi \downarrow \\ \overline M_{0,n} \end{array}$$
where $ev(f) = (f(p_1),...,f(p_n))$ and $\pi(f)$ is the ``stabilization'' of 
$(C;p_1,...,p_n)$. Of course, $\pi$ is only defined when
$n\ge 3$. Finally, there is a diagram:
$$\begin{array}{ccc}\overline M_{0,n+1}(X,\beta) & \stackrel {e}\rightarrow & X\\
f\downarrow
\\
\overline M_{0,n}(X,\beta)\end{array}$$
where $f$ forgets the marked point $p_{n+1}$ and stabilizes, and
$e(f) = f(p_{n+1})$. This exhibits an isomorphism between 
$\overline M_{0,n+1}(X,\beta)$ and the ``universal curve'' ${\cal C}$.

\medskip

When $X$ is a ``nice enough'' (see the theorem), then we define:
$$c^\beta_{0,n}(X) = ev_*(1)$$

Otherwise, the ``right'' definition for 
$c^\beta_{0,n}(X)$ has been found by Li-Tian and Behrend-Fantechi
via the ``virtual'' fundamental class on $\overline M_{0,n}(X,\beta)$
\cite{BehF,LiT}.
Here I will ignore these subtleties, referring to these classes only as $ev_*(1)$
and using the properties listed in \cite{BehM} wherever convenient.

\medskip

A couple of key relations among Gromov-Witten invariants follow:

\medskip

\nt {\bf Basic Axioms:} 
(a) {\bf (Easy String)} $\langle \gamma_1,...,\gamma_n,1\rangle_\beta = 0$
unless $n=2,\beta=0$.

\medskip

(b) {\bf (Easy Divisor)} If $D\in V$ is the class of a divisor on $X$, then 

$$\langle
\gamma_1,...,\gamma_n,D\rangle_\beta = 
(D\cap \beta)\langle \gamma_1,...,\gamma_n\rangle_\beta \ \ \mbox{unless $n=2,\beta
= 0$}$$ 

\medskip

{\bf Proof:} Apply the projection formula to the commuting diagram:
$$\begin{array}{ccc} \overline M_{0,n+1}(X,\beta) & 
\stackrel {ev}\rightarrow & X^{n+1} \\
f \downarrow & & \downarrow \\ \overline M_{0,n}(X,\beta) & \stackrel{ev}
\rightarrow & X^n
\end{array}$$
together with the fact that $f_*1 = 0$ and $f_*e^*D = (D\cap \beta)\cdot 1$.

\medskip

The proof of WDVV requires the additional map to $\overline M_{0,n}$.
Given a subset $S\subseteq \{1,...,n\}$ of cardinality $m$ and
$\alpha\in C_{eff}$, there is a gluing morphism:
$$\delta_{S,\alpha}: \overline M_{0,m+1}(X,\alpha) \times_X \overline
M_{0,n-m+1}(X,\beta -
\alpha) \rightarrow \overline M_{0,n}(X,\beta)$$
sending a pair of maps $f:(C;p_1,...,p_m,p) \rightarrow X$ and
$g:(B;q_1,...,q_{n-m},q) \rightarrow X$ satisfying $f(p) = g(q)$ to
the single map $h:(C\cup_{p=q} B; p_i,q_j) \rightarrow X$. 
To give the ordering of the points on $C\cup B$, apply the
permutation $\sigma_S$ to $\{p_i,q_j\}$ given by:
$$\sigma_S(i) = s_i,\sigma_S(m+i) = s^c_i\ \mbox{where}\ S =
\{s_1<...<s_m\}, S^c = \{s^c_1<...<s^c_{n-m}\}$$

The gluing maps have two important properties:

\medskip

$\bullet$ $ev_*{\delta_{S,\alpha}}_*(1) =
\sigma_S(c^\alpha_{0,m+1}\rfloor c^{\beta -\alpha}_{0,n-m+1})$
where $\sigma_S$
permutes the factors of $V^{\otimes n}$.

\medskip

$\bullet$ $\pi^*{\delta_S}_*(1) = \sum_\alpha {\delta_{S,\alpha}}_*(1)$
for $\delta_S:\overline M_{0,m+1}\times \overline
M_{0,n-m+1}\rightarrow \overline M_{0,n}$.

\medskip

In addition, one easily checks that:

\medskip

$\bullet$ $cr^*{\delta_{\{ij\}}}_*(1) = \sum_{S}
{\delta_{S \cup \{n+i,n+j\}}}_*(1)$ for $i,j\le 4$ and the
forgetful cross-ratio map
$cr:\overline  M_{0,n+4}\rightarrow \overline M_{0,4}$
But ${\delta_{\{ij\}}}_*(1)$ is the class of a point on $\overline M_{0,4} = {\bf
P}^1$, so is independent of $i,j$, and
all together, we get:
$$\sum_m\sum_{|S| = m}\sum_\alpha \sigma_{S\cup
\{n+1,n+2\}}(c^\alpha_{0,m+3}\rfloor 
c^{\beta - \alpha}_{0,n-m+3}) = 
ev_*\pi^*cr^*{\delta_{\{12\}}}_*(1) =$$
$$= ev_*\pi^*cr^*{\delta_{\{13\}}}_*(1) = 
\sum_m\sum_{|S| = m}\sum_\alpha \sigma_{S\cup
\{n+1,n+3\}}(c^\alpha_{0,m+3}\rfloor 
c_{0,n-m+3}^{\beta - \alpha})$$
which are precisely the coefficients of $q^\beta$ in 
$c_0(X)\rfloor c_0(X)$
and the pull-back of $c_0(X)\rfloor c_0(X)$ under the $(2\ 3)$
transposition,  respectively.

\medskip

\nt {\bf Remark:} Checking these properties of the gluing maps
with virtual classes is the major 
technical point in the work of Behrend-Fantechi and Li-Tian.

\bigskip

\nt {\bf 2. Kontsevich-Manin Reconstruction.} This is generally understood to
be the following consequence of the WDVV relations (and the basic axioms):

\medskip

\nt {\bf Theorem:} (Kontsevich-Manin \cite{KonM}) If $V$ is generated by
divisor classes, then
$c_0(X)$ can be explicitly reconstructed from the classes
$c^\beta_{0,2}(X)$ and $c^0_{0,3}(X) = \Delta_3$.  

\medskip

{\bf Proof:} Choose generating divisor classes $D_1,...,D_k$, 
and choose a basis $\{D_J = \prod D_i^{j_i}\} \in V$ consisting of monomials in the
divisors $D_j$.
For $n \ge 0$,$\beta > 0$ and $|K| \ge 2$, we use WDVV to 
solve for Gromov-Witten invariants of the form: 
$$\langle \gamma_1,...,\gamma_{n+2},D_K \rangle_\beta$$
in terms of invariants involving either:

\medskip

$(i)$ smaller $\beta$ (and either the same or smaller $n$) or 

\medskip

$(ii)$ the same $\beta$ and $n$ but
smaller $|K|$. 

\medskip

Since the basic axioms allow us to eliminate $D_K$ in case $|K| =
0$ or $1$, the theorem then follows by induction.

\medskip

Let $(g^{IJ})$ be the inverse of the intersection matrix defined by:
$$g_{IJ} := \int_X D_I \wedge D_J$$
and choose $D_L,D_M$ so that
$D_L \wedge D_M = D_K$ with $|L|,|M| < |K|$. Then the coefficient of 
$q^\beta$ in the expression:
$$\int_{X^{n+4}} c_0(X)\rfloor c_0(X) \wedge 
\gamma_1\otimes ... \otimes \gamma_{n+2}\otimes D_L\otimes D_M$$
is, explicitly,
$$(1)\ \sum_S\sum_{I,J}\langle
\gamma_{s_1},...,\gamma_{s_m},\gamma_{n+1},\gamma_{n+2},D_I\rangle_{\alpha}
g^{IJ} \langle D_J,\gamma_{s^c_1},...,\gamma_{s^c_{n-m}},
D_L,D_M\rangle_{\beta - \alpha}$$ 
whereas after the $(2\ 3)$ permutation, we obtain:
$$(2)\ \sum_S\sum_{I,J}\langle
\gamma_{s_1},...,\gamma_{s_m},\gamma_{n+1},D_L,D_I\rangle_{\alpha}
g^{IJ} \langle D_J,\gamma_{s^c_1},...,\gamma_{s^c_{n-m}},
\gamma_{n+2},D_M\rangle_{\beta - \alpha}$$ 
Only the terms $\alpha = 0, S = \emptyset$ or $\alpha = \beta$,
$S = \{1,...,n\}$ do not satisfy $(i)$.
But these terms can be computed using the fact
that $c^0_{0,3}(X) = \Delta_3$. They are:

\medskip

$(1)\ \langle \gamma_{n+1}\wedge
\gamma_{n+2},\gamma_1,...,\gamma_n,D_L,D_M\rangle_\beta + \langle
\gamma_1,...,\gamma_{n+2},D_K\rangle_\beta$ and  

\medskip

$(2)\ \langle \gamma_{n+1}\wedge
D_L,\gamma_1,...,\gamma_n,\gamma_{n+2},D_M\rangle_\beta + \langle
\gamma_1,...,\gamma_{n+1},D_L,\gamma_{n+2} \wedge
D_M\rangle_\beta$

\medskip

\nt The theorem follows since we are solving for the second of
the terms in $(1)$ and the other three of these ``special'' terms satisfy
$(ii)$.

\bigskip

\nt {\bf Examples:} (a) For $X = {\bf P}^3$ and $\beta =1$, there
are only the special terms, and taking $\gamma_1 = H^2,\gamma_2 =
H^2, \gamma_3 = H^2,D_L=H,D_M=H$ gives:
$$\langle H^4,H^2,H,H\rangle_1 + \langle H^2,H^2,H^2,H^2\rangle_1
= \langle H^3,H^2,H^2,H\rangle_1 + \langle
H^2,H^2,H,H^3\rangle_1$$
from which one concludes using the divisor axiom that:
$$\langle H^2,H^2,H^2,H^2\rangle_1 = 2\langle
H^3,H^2,H^2\rangle_1$$
and letting $\gamma_1 = H^3,\gamma_2 = H^2,D_L=H,D_M=H$ gives:
$$\langle H^5,H,H\rangle_1 + \langle H^3,H^2,H^2\rangle_1 = 
\langle H^4,H^2,H\rangle_1 + \langle H^3,H,H^3\rangle_1$$
which gives $\langle H^3,H^2,H^2\rangle_1 = \langle
H^3,H^3\rangle_1 = 1$ and $\langle H^2,H^2,H^2,H^2\rangle_1 = 2$.

\medskip

(b) For a more interesting example, let
let $n_d$ be the
number of rational plane curves of degree $d$
through $3d-1$ general points of ${\bf P}^2$. That is,
$$n_d = \langle H^2,....,H^2\rangle_d$$

Let $\gamma_1,...,\gamma_{3d-2} = H^2$ and $D_L=D_M=H$. Then:
$$(1)\ \ 0 + n_d + 
\sum_{0<e<d}\sum_{|S| = 3e-3}\langle H^2,...,H^2,H^2,H^2,H\rangle_e
\langle H,H^2,...,H^2,H,H\rangle_{d-e}$$ 
$$(2)\ \ 0 + 0 + \sum_{0<e<d}\sum_{|S|=3e-2}
\langle H^2,...,H^2,H^2,H,H\rangle_e
\langle H,H^2,...,H^2,H^2,H\rangle_{d-e}$$
and this, together with the divisor equation, gives the recursive
formula:

$$n_d = 
\sum_{e=1}^{d-1}n_en_{d-e}\left(e^2(d-e)^2\left({3d-4} \atop
{3e-2}\right) - e(d-e)^3\left({3d-4}\atop
{3e-3}\right)\right)$$ 
from which (together with the a priori $n_1 = 1$) one computes:
$$n_2 = 1, n_3 = 12, n_4 = 620 \ (\mbox{as in the introduction})$$

\nt {\bf Remarks:} The reconstruction theorem applies in a more general
situation. If we can write $V = W + W^{\perp}$ so that $W$ is generated
by divisor classes and every Gromov-Witten invariant 
invariant of the form $\langle
\gamma_1,...,\gamma_n,\lambda\rangle_\beta$ is zero when
$\gamma_1,...,\gamma_n\in W$ and
$\lambda \in W^{\perp}$, then every Gromov-Witten invariant of the form
$\langle \gamma_1,...,\gamma_n\rangle_\beta$ (for all $\gamma_i\in W$)
can be reconstructed from the $n=2$ invariants.
The same proof gives this result, as well.

\medskip

This applies, for example, to the case of a hypersurface 
$X\subset {\bf P}^n$. Which begs the question: What are the
$c^{\beta}_{0,2}(X)$ in that case? Holger Kley and I found a way to
reconstruct these, too, using other hidden symmetries.

\medskip

\nt {\bf 3. A Symmetry of $J$-functions.} The universal curve:
$$\begin{array}{ccc} {\cal C} \\ f\downarrow \\ \overline M_{0,n}(X,\beta)
\end{array}$$
with universal sections (the marked points) $\rho_i:\overline M_{0,n}(X,\beta)
\rightarrow {\cal C}$ determine the
``Morita classes''
$\psi_i := c_1(N^*_{\rho_i})$, where $N^*_{\rho_i}$ is the conormal bundle
of the section. The following strange-looking definition is 
inspired by mirror symmetry
(as interpreted by Givental \cite{Giv}, Lian-Liu-Yau \cite{LLY} and
others):  
$$J^\beta_{0,n}(X) := ev_*\left(\frac 1{\prod_{i=1}^n
t_i(t_i-\psi_i)}\right)
\in Sym^n(V[t^{-1}])$$
where $t$ (or rather $t_i$) is a variable, and the denominator is formally inverted.
Note that the coefficient of
$\prod t_i^{-2}$ is $c^\beta_{0,n}(X)$.

\medskip

These seem to have better properties than the
$c^\beta_{0,n}(X)$.
For example, the $n=1$ (mirror conjecture) case has some very interesting 
``functorial'' properties and is often amenable to computations.

\medskip

\nt {\bf Theorem:} (Bertram/Behrend \cite{Ber,Beh2}) If $\beta\in
C_{eff}(X\times Y)$, let
$\beta_1\in C_{eff}(X)$ and $\beta_2\in C_{eff}(Y)$ be the two
projections. Then:
$$J^{\beta}_{0,1}(X\times Y) = J^{\beta_1}_{0,1}(X)\otimes_{{\bf
Q}[t^{-1}]} J^{\beta_2}_{0,1}(Y)$$
(there is a product formula for $c^\beta_{0,2}(X)$ classes, but it is
{\bf much} more complex)

\medskip

\nt {\bf Theorem:} (Givental \cite{Giv})  Let $H$ be the hyperplane
class in
$H^2({\bf P}^n,{\bf Z})$. Then:
$$J^d_{0,1}({\bf P}^n) = \frac 1{\prod_{k=1}^d(H+kt)^{n+1}} $$
If $X\subset {\bf P}^n$ is a complete intersection of type
$(l_1,...,l_r)$ of dimension
$\ge 3$ with $l_1+...+l_r < n$, then:
$$J^d_{0,1}(X) = \frac
{\prod_{i=1}^r
\prod_{k=1}^{dl_i}(l_iH+kt)}{\prod_{k=1}^d(H+kt)^{n+1}}$$
and if$l_1+...+l_r = n$ or $n+1$, then $J^d_{0,1}(X)$ is 
computed from the $J^d_{0,1}({\bf P}^n)$ by an  explicit ``mirror
transformation'' (see \cite{Ber2}).

\medskip

This generalizes to complete intersections in (Fano) toric varieties,
but it is even more general, as pointed out first by Kim and proved
by Lee \cite{Lee}:

\medskip

\nt {\bf Quantum Lefschetz Hyperplane Theorem:} If $X \subset Y$ is a Fano
or Calabi-Yau very ample hypersurface of dimension at least $3$, then the
$J^\beta_{0,1}(X)$ are explicitly determined
by the $J^{\beta}_{0,1}(Y)$.

\medskip

The following computation:
$$J^0_{0,n}(X) = \Delta_n\prod t_i^{-2}\left(\sum t_i^{-1}\right)^{n-3}$$
follows from the well-known intersection numbers on $\overline M_{0,n}$:
$$\int_{\overline M_{0,n}} \psi_1^{a_1}\wedge ...\wedge \psi_n^{a_n} = 
\frac {(n-3)!}{a_1!...a_n!}$$
whenever $n \ge 3$ and $\sum a_i = n-3$. We argue by analogy that:
$$J^0_{0,1}(X) = 1 \ \ \mbox{and}\ \ \ J^0_{0,2}(X) = \frac
{\Delta_2}{t_1t_2(t_1+t_2)}$$
are good definitions (though the second one does not belong to
$Sym^2(V[t^{-1}])$).

With these definitions, we put together the generating function:
$$J_0(X) := \sum_{\beta}\sum_{n \ge 1} J^\beta_{0,n}(X)q^\beta \in
Sym^*(V[t^{-1}])[[q]]$$ (with a little extra room for $J^0_{0,2}(X)$), and
by analogy with WDVV, we have:

\medskip

\nt {\bf Theorem 1:} (Bertram-Kley \cite{BerK}) The $J_0(X)$
invariants have the symmetry:
$$J_0(X)\rfloor J_0(X) = 0$$
defined as follows. Let $U = V[t^{-1}]$ and extend Poincar\'e duality
to:
$$U[[q]] \otimes U[[q]] \rightarrow {\bf
Q}[t^{-1}][[q]];
\gamma_1 t^{-k}q^\alpha
\otimes
\gamma_2 t^{-l}q^\beta \mapsto \left(\int_X \gamma_1 \wedge \gamma_2\right)
t^{-k}(-t)^{-l}q^{\alpha+\beta}$$ 
Partially desymmetrize, writing $J_0(X)$ as an element of
$Sym^*(U)\otimes U[[q]]$, and:
$$J_0(X)\rfloor J_0(X) = 0 \in Sym^*(U)[t^{-1}][[q]]$$
is a collection of identities on the coefficients of the $t^{-k}$.

\medskip

\nt {\bf The Extra Term:} Contracting with $J^0_{0,2}(X)$ does not fit
the pattern above, at least not in the obvious way. Instead, we
extend the Poincar\'e duality to:
$$\gamma_1t^{-k}\otimes \frac{\gamma_2t^{-l}}{t_i+t} \mapsto
\left(\int_X\gamma_1\wedge \gamma_2\right)\frac{t^{-k} -
t_i^{-k}}{t_i-t}(-t)^{-l}$$ 
so that $J^0_{0,2}(X)$ acts as a difference operator, and similarly,
$$\frac{\gamma_1t^{-k}}{t_i+t}\otimes \gamma_2t^{-l} \mapsto
\left(\int_X\gamma_1\wedge
\gamma_2\right)\frac{(-t)^{-l} - t_i^{-l}}{t_i+t}t^{-k}$$ 
As an example, we use Theorem 1 to prove the:

\medskip

\nt {\bf String Equation:} Define the 
generalized Gromov-Witten invariants by:
$$\langle
\psi^{a_1}(\gamma_1),...,\psi^{a_n}(\gamma_n)\rangle_\beta :=
\int_{X^n}ev_*(\psi_1^{a_1}...\psi_n^{a_n}) \wedge \gamma_1 \otimes
...\otimes \gamma_n$$ if each $a_i \ge 0$, and zero if some $a_i <
0$. Then unless $n=1$ or $2$ and $\beta = 0$,
$$\langle
\psi^{a_1}(\gamma_1),...,\psi^{a_n}(\gamma_n),1\rangle_\beta =
\sum_{i=1}^n\langle
\psi^{a_1}(\gamma_1),...,
\psi^{a_i-1}(\gamma_i),...,\psi^{a_n}(\gamma_n)\rangle_\beta$$

{\bf Proof:} The $U^{\otimes n}q^\beta$ coefficient of $J_0(X)\rfloor
J_0(X)$ is:
$$\sum_\alpha \sum_m \sum_{|S|=m} \sigma_S\left(J^\alpha_{0,m+1}(X)
\rfloor J^{\beta - \alpha}_{0,n-m+1}(X)\right)$$ 
and we say that a summand involving
$J^0_{0,1}(X)$ or $J^0_{0,2}(X)$ is special. 
Let $\{e_i\}$ be a basis of $V$ with intersection
matrix $g_{ij}$, and consider the 
$q^\beta$ coefficient of:
$$\int_{X^n} J_0(X)\rfloor J_0(X) \wedge \gamma_1\otimes ...\otimes
\gamma_n$$
The non-special summands produce terms of the form:
$$\sigma_S\langle
\frac{\gamma_1}{t_1(t_1-\psi)},...,\frac{e_i}{t(t-\psi)}\rangle_\alpha
g^{ij}\langle
\frac{e_j}{-t(-t-\psi)},...,\frac{\gamma_n}{-t_n(-t_n-\psi)}\rangle_{\beta
- \alpha}$$
where we formally set
$\frac{\gamma}{t(t-\psi)} = t^{-2}\gamma + t^{-3}\psi(\gamma) +
t^{-4}\psi^2(\gamma) + ...$

On the other hand, the special summands produce:
$$(m=n,\alpha = \beta)\ \ \ \  \langle
\frac{\gamma_1}{t_1(t_1-\psi)},...,\frac{\gamma_n}{t_n(t_n-\psi)},
\frac 1{t(t-\psi)}\rangle_\beta\ \ \ \mbox{and}$$ 
$$(m=n-1,\alpha = \beta)\ \ \ \ 
\sum_{i=1}^n\langle
\frac{\gamma_1}{t_1(t_1-\psi)},...,
f_{t,t_i}(\gamma_i),...,\frac{\gamma_n}{t_n(t_n-\psi)}\rangle_\beta$$
where 
$$f_{t,t_i} = \frac{(t(t-\psi))^{-1} -
(t_i(t_i-\psi))^{-1}}{t_i(-t)(t_i - t)}$$ 
and the $\alpha = 0$ terms, which are the same functions, but of
$-t$ instead of $t$.  

Setting the coefficient of $t^{-2}$ to zero (using Theorem 1) gives:
$$\langle
\frac{\gamma_1}{t_1(t_1-\psi)},..,\frac{\gamma_n}{t_n(t_n-\psi)},
1\rangle_\beta =  
\sum_{i=1}^n\langle
\frac{\gamma_1}{t_1(t_1-\psi)},...,
\frac{\gamma_i}{t_i^2(t_i-\psi)},...,\frac{\gamma_n}{t_n(t_n-\psi)}\rangle_\beta$$
which is the string equation!

\medskip

\nt {\bf 4. More Symmetries of $J$-Functions.} Theorem $1$ is a consequence
of a family of symmetries of $J$-functions, which will occupy the rest
of this paper.

\medskip

\nt {\bf Definition:} Let $J'_0(X) \in Sym^*U\otimes U^{\otimes 2}[[q]]$
be the desymmetrization of $J_0(X)$.

\medskip

(this notation distinguishes it from $J_0(X) \in
Sym^*U\otimes U[[q]]$) 

\medskip

\nt {\bf Theorem 2:} (Bertram-Kley \cite{BerK})  
$$J'_0(X)\rfloor J_0(X) = 0 \in Sym^*U\otimes U[t^{-1}][[q]]$$
There are more symmetries, but first I want to give a:

\medskip

\nt {\bf Stringy Corollary:} Every (generalized) Gromov-Witten invariant
of the form:
$$\langle \psi^{a_0}(\gamma_0),
\psi^{a_1}(\gamma_1),...,\psi^{a_n}(\gamma_n),\psi^{a_{n+1}}(1)\rangle_\beta$$
can be explicitly expressed in terms of ``simpler'' invariants.

\medskip

Consider the $q^\beta$ coefficient of:
$$\int_{X^{n+1}}J_0'(X)\rfloor J_0(X) \wedge \gamma_0\otimes ...\otimes
\gamma_n$$
(where we index the distinguished factor of $Sym^*U\otimes U$ with the zero
subscript). 

\medskip

This consists of special terms:
$$\langle
\frac{\gamma_0}{t_0(t_0-\psi)},...,\frac{\gamma_n}{t_n(t_n-\psi)},\frac
1{t(t-\psi)}\rangle_\beta + 
\langle
f_{-t,t_0}(\gamma_0),\frac{\gamma_1}{t_1(t_1-\psi)}...,
\frac{\gamma_n}{t_n(t_n-\psi)}\rangle_\beta$$
$$+ \sum_{i=1}^n \langle
\frac{\gamma_0}{t_0(t_0-\psi)},...,f_{t,t_i}(\gamma_i),\frac{\gamma_n}{t_n(t_n-\psi)},\frac
1{t(t-\psi)}\rangle_\beta$$
where $f_{t,t_i}$ is the divided difference of $(t_i(t_i-\psi))^{-1}$
defined earlier.

\medskip

In addition, there are the non-special terms, of the form:
$$\sigma_{\{0\}\cup S}
\langle
\frac{\gamma_0}{t_0(t_0-\psi)},...,\frac
{e_i}{t(t-\psi)}\rangle_\alpha g^{ij}
\langle \frac{e_j}{-t(-t-\psi)},
,...,\frac{\gamma_n}{-t_n(-t_n-\psi)}
\rangle_{\beta - \alpha}$$

But Theorem 2 says that the sum of all the terms is zero. So since the first
special term generates all the Gromov-Witten invariants of the corollary, and
all the other terms are simpler, the corollary is proved.

\medskip

\nt {\bf Examples:} As before, we get the string equation from the
coefficient of $t^{-2}$.
From the coefficient of $t^{-3}$, we obtain the {\bf Dilaton Equation:}
$$\langle \psi^{a_0}(\gamma_0),
...,\psi^{a_n}(\gamma_n),\psi(1)\rangle_\beta = (n-1)
\langle
\psi^{a_0}(\gamma_0),...,\psi^{a_n}(\gamma_n)\rangle_\beta$$

Fix divisor classes $H_1,...,H_n\in V$ and a polynomial $p(x_1,...,x_n)$.
For each $\beta \in C_{eff}$, let $b_i = H_i \cap
\beta$, and  define:  
$$p(H-\beta t) = p(H_1 - b_1t,...,H_n - b_nt)$$
Then we have the following:

\medskip

\nt {\bf Theorem 3:} (Bertram-Kley \cite{BerK}) For any polynomial
$p$ as above:
$$J'_0(X)\rfloor_{p} J_0(X) \in
Sym^*(U)\otimes U[[q]]\otimes t^{-1}{\bf Q}[t]$$
where ``$\rfloor_{p}$'' is the contraction defined using the pairing:
$$\gamma_1 t^{-k}q^\alpha \otimes \gamma_2 t^{-l}q^\beta
\mapsto \left(\int_X \gamma_1 \wedge \gamma_2\wedge p(H-\beta t)\right)
t^{-k}(-t)^{-l}
\in {\bf Q}[t,t^{-1}]q^{\alpha+\beta}$$

\nt {\bf Remark:} Theorem 3 $\Rightarrow$ Theorem 2 ($\Rightarrow$ Theorem
1) taking $p = 1$.

\medskip

\nt {\bf Example:} The {\bf Divisor Equation} is obtained from $p(x)=x$.
The special terms of:
$$\int_{X^{n+1}} J'_0(X)\rfloor_xJ_0(X) \wedge \gamma_0\otimes
...\otimes \gamma_n$$
are:
$$\langle
\frac{\gamma_0}{t_0(t_0-\psi)},...,\frac
H{t(t-\psi)}\rangle_\beta + 
\langle
f_{-t,t_0}(\gamma_0\wedge (H-bt)),\frac{\gamma_1}{t_1(t_1-\psi)}...,
\frac{\gamma_n}{t_n(t_n-\psi)}\rangle_\beta$$
$$+ \sum_{i=1}^n \langle
\frac{\gamma_0}{t_0(t_0-\psi)},...,f_{t,t_i}(\gamma_i\wedge
H),\frac{\gamma_n}{t_n(t_n-\psi)},\frac 1{t(t-\psi)}\rangle_\beta$$
and  
non-special terms do not contribute to the coefficient of $t^{-2}$, which
gives:
$$\langle \psi^{a_0}(\gamma_0),
...,\psi^{a_n}(\gamma_n),H\rangle_\beta = 
b\langle
\psi^{a_0}(\gamma_0),...,\psi^{a_n}(\gamma_n)\rangle_\beta$$
$$+ \sum_{i=0}^n\langle
\psi^{a_0}(\gamma_0),...,\psi^{a_i-1}(\gamma_i\wedge H),...,
\psi^{a_n}(\gamma_n)\rangle_\beta
$$

\nt {\bf Reconstruction Corollary:} If $V$ is generated by divisor 
classes, then $J_0(X)$ can be explicitly reconstructed from the
$J^\beta_{0,1}(X)$ classes.

\medskip

{\bf Proof:} When we pair $J'_0(X)\rfloor_p J_0(X)$ with $\gamma_0\otimes
...\otimes \gamma_n$, then the first special term is of the form:
$$\langle
\frac{\gamma_0}{t_0(t_0-\psi)},...,\frac{\gamma_n}{t_n(t_n-\psi)},\frac
{p(H)}{t(t-\psi)}\rangle_\beta$$
and all the other terms are ``simpler.'' Since this only involves
powers $t^{-2},t^{-3},...$, it is inductively determined by
the $J^\beta_{0,1}(X)$ using Theorem 3 and the fact that $V$ is
generated by divisor classes.  

\medskip

\nt {\bf Remark:} As with Kontsevich-Manin reconstruction, this 
holds in greater generality. Let $W \subset V$ be the subalgebra generated
by divisor classes. If
the coefficients of $J^\beta_{0,1}(X)$ all belong to $W$, then
all generalized invariants involving only classes coming from $W$
(i.e. the projection of $J_0(X)$ to $Sym^*(W[t^{-1}])[[q]]$) can be
explicitly reconstructed from the $J^\beta_{0,1}(X)$.

\medskip

If we only consider the $\prod t_i^{-2}$ terms, then by Theorem
3, the coefficients of $t^{-2},t^{-3},...$ vanish in 
each of the following
expressions:
$$\langle \gamma_0,...,\gamma_n,\frac{p(H)}{t(t-\psi)}\rangle_\beta
+ \sum_{i=1}^n \langle \gamma_0,..,\frac{\gamma_i\wedge p(H)}
{-t^2(t-\psi)},..,\gamma_n\rangle_\beta + \langle
\frac{\gamma_0\wedge
p(H-bt)}{-t^2(-t-\psi)},...,\gamma_n\rangle_\beta$$
$$+ \sum_\alpha \sum_m\sum_{|S|=m}\sigma_{\{0\}\cup S}\langle
\gamma_0,...,\gamma_m,\frac{e_i}{t(t-\psi)}\rangle_{\beta - \alpha}
g^{ij}\langle \frac{e_j\wedge
p(H-at)}{-t(-t-\psi)},\gamma_{m+1},...,\gamma_n\rangle_\alpha$$
and this leads to very efficient algorithms for reconstructing
ordinary invariants, provided $J^\beta_{0,1}(X)$ are known and
the conditions of the remark above are satisfied. Note that
in order to reconstruct ordinary invariants, we can stay within
the realm of invariants with at most one ``gravitational descendant''
$\psi^a$.

\medskip

\nt {\bf Example:} Another look at rational plane curves.
Let $X = {\bf P}^2$ and:
$$n_d^{(a)} := \langle H^2,...,H^2,\psi^a(H^2)\rangle_d\ \ \mbox{
(the total number of terms is $3d-1-a$)}$$

Then it follows from our relations that:
$$\begin{array}{lcl} n_d & = & d^2n_d^{(1)} - \sum_{e=1}^{d-1} \left(
{3d-3}\atop {3e-1}\right) (d-e)e^3n_{d-e}n_e\\
n_d^{(1)} & = & d^2n_d^{(2)} -
\sum_{e=1}^{d-1}\left({3d-4}\atop{3e-1}\right)
(d-e)e^3n_{d-e}^{(1)}n_e -
\sum_{e=1}^{d-1}\left({3d-4}\atop{3e-2}\right) e^2n_{d-e}n_e\\
n_d^{(2)} & = & d^2n_d^{(3)} -
\sum_{e=1}^{d-1}\left({3d-5}\atop{3e-1}\right)
(d-e)e^3n_{d-e}^{(2)}n_e -
\sum_{e=1}^{d-1}\left({3d-5}\atop{3e-2}\right) e^2n_{d-e}^{(1)}n_e\\
\ \ \vdots \\
n_d^{(3d-3)} & = & d^2n_d^{(3d-2)} = \frac{1}{d(d-1!)^3} \
\mbox{(from Givental)}
\end{array}$$
In case $d\le 3$, this gives: $n_1 = n^{(1)}_1 = 1; n_2 =
n^{(1)}_2 = n^{(2)}_2 = 1, n^{(3)}_2 = \frac 12, n^{(4)}_2 = \frac
18;$
$$n_3 = 12, n_3^{(1)} = 10,
n_3^{(2)} = 7, n_3^{(3)} = 3, n_3^{(4)} = 1, n_3^{(5)} = \frac 14,
n_3^{(6)} = \frac 1{24},n_3^{(7)} = \frac 1{216}$$
It is amusing to note that
$d=1$ ``proves'' there is one line through $2$ points. But
seriously, the ``mirror data'' $J^\beta_{0,1}(X)$ does appear in
general to be more basic  and a better starting point for inductions
than the two-point invariants $c^\beta_{0,2}(X)$.
On the other hand, if one is only interested in ordinary invariants,
it is probably more efficient to compute the $c^\beta_{0,2}(X)$
first and then to apply Kontsevich-Manin reconstruction. This is easy to
implement as an algorithm when $X$ is a Fano hypersurface  
in ${\bf P}^n$ (see the appendix in \cite{BerK}).

\newpage

\nt {\bf 5. Localization.} The proof of Theorem 3 relies on the localization
theorem of Atiyah-Bott, as adapted by Graber and Pandharipande
\cite{GraP} to apply to virtual classes on the moduli stacks of stable maps.
The idea of using localization to compute Gromov-Witten invariants was
introduced by  Kontsevich \cite{Kon}, but the key insight, described below, 
belongs to Givental \cite{Giv}. 

\medskip

Recall that in order to prove WDVV, one considered the map:
$$\begin{array}{c}\overline M_{0,n}(X,\beta) \\ \downarrow \\
\overline M_{0,n}\end{array}$$
and relations in $H^*(\overline M_{0,n},{\bf Q})$ gave
relations among Gromov-Witten invariants.

\medskip

Here, the key insight is to consider the embedding:
$$i:\overline M_{0,n}(X,\beta) \hookrightarrow \overline M_{0,0}(X\times 
({\bf P}^1)^n,(\beta,1^n))$$
replacing each marked point with a curve mapping isomorphically to ${\bf
P}^1$, with the node mapping to $0$ (see \cite{Ber2}). In this way,
$\overline M_{0,n}(X,\beta)$ is one of the (many!) connected components of
the locus of fixed points for the natural action of the torus $T = ({\bf
C}^*)^n$ on the ``graph space''
$\overline M_{0,0}(X\times ({\bf P}^1)^n,(\beta,1^n))$. 

The localization 
theorem says that an equivariant cohomology (or Chern) class $c$ on the graph 
space can be recovered from the total fixed-point locus $i:F \hookrightarrow 
\overline M_{0,0}(X\times ({\bf P}^1)^n,(\beta,1^n))$
via the localization formula:
$$c \equiv i_*\frac{i^*c}{e_T(F)} \ \ \mbox{mod torsion}$$
($e_T(F)$ is the equivariant Euler class) and one can compute fairly easily
that:
$$e_T(\overline M_{0,n}(X,\beta)) = \prod_{i=1}^n t_i(t_i - \psi_i)$$
where $H^*(BT,{\bf Q}) = {\bf Q}[t_1,...,t_n]$. 

One could conceivably use this idea to compute Gromov-Witten invariants, since
all the components of $F$ are fiber products of lower-degree stable map
spaces, and the Euler classes are always expressible in terms of 
the $t_i$ and $\psi_i$. But the
combinatorial problem of enumerating the components of $F$ is a nightmare
(see \cite{Kon}). 
Instead, Kley and I use forgetful maps among
these graph spaces:
$$f:\overline M_{0,0}(X\times ({\bf P}^1)^{n+1},(\beta,1^{n+1})) \rightarrow
\overline M_{0,0}(X\times ({\bf P}^1)^{n},(\beta,1^{n}))$$
and
$$g:\overline M_{0,0}(X\times ({\bf P}^1)^{n+1},(\beta,1^{n+1})) \rightarrow
\overline M_{0,0}({\bf P}^1\times {\bf P}^1,(1,1)) \cong {\bf P}^3$$
to cobble together a birational $T$-equivariant map:
$$\Phi: M_{0,0}(X\times ({\bf P}^1)^{n+1},(\beta,1^{n+1})) \rightarrow
\overline M_{0,0}(X\times ({\bf P}^1)^{n},(\beta,1^{n}))\times {\bf P}^3$$

The advantage of this map is that the connected components of the total fixed
locus $i':F' \hookrightarrow M_{0,0}(X\times ({\bf
P}^1)^{n+1},(\beta,1^{n+1}))$ lying in the preimage of 
$j:M_{0,n}(X,\beta)
\hookrightarrow M_{0,0}(X\times ({\bf P}^1)^{n},(\beta,1^{n})) \times {\bf
P}^3$ are not hard to enumerate. And from the more refined localization
formula:
$$\Phi_*(\frac {i'^*c}{e_T(F')}) = \frac{j^*\Phi_*c}{e_T(\overline
M_{0,n}(X,\beta)} = 
\frac{j^*\Phi_*c}{t_1t_{n+1}(t_1+t_{n+1})\prod t_i(t_i-\psi_i)}$$
we were able to deduce Theorem 3, by further pushing forward to 
$X^n$ via $ev$ and choosing suitable equivariant Chern 
classes $c$ on the graph space. 

\medskip

This approach does have its limitations, however. It does not suggest
a ``reason'' for the particularly nice organization of the relations in Theorem
3. In particular, the same approach yields many relations among $J$-functions
for higher genus curves, which ought to be organized and analyzed in the
context of the many conjectures about generating functions
made up of higher-genus Gromov-Witten invariants. Perhaps 
the physicists can come to our rescue again.

\nt email: bertram@math.utah.edu


\begin{thebibliography}{10}

\bibitem{Beh}
K. Behrend, {\it Gromov-Witten invariants in algebraic 
geometry}, Invent Math {\bf 127} (1997) 601--617.

\bibitem{Beh2} K. Behrend, {\it The product formula for Gromov-Witten
invariants}, J Alg Geom {\bf 8} (1999) 529-541.

\bibitem{BehF}
K. Behrend and B. Fantechi, {\it The intrinsic normal cone},
Invent Math {\bf 128} (1997) 45--88.

\bibitem{BehM}
K. Behrend and Y. Manin, {\it Stacks of stable maps and Gromov-Witten
invariants}, Duke Math J {\bf 85} (1996) 1--60.

\bibitem{Ber} A. Bertram,
{\it Some applications of localization to enumerative problems},
Michigan Math J {\bf 48} (2000) 65-75.

\bibitem{Ber2}
A. Bertram, {\it Another way to enumerate rational curves with torus actions},
Invent Math {\bf 142} (2000) 487-512.

\bibitem{BerK} A. Bertram and H. Kley, {\it New recursions for
genus-zero Gromov-Witten invariants}, math.AG/0007082.

\bibitem{Can} P. Candelas, X. de la Ossa, P. Green and L. Parkes,
{\it A pair of Calabi-Yau manifolds as an exactly soluble superconformal
field theory}, Nuclear Phys. B {\bf 359} (1991) 21-74.


\bibitem{Ell} G. Ellingsrud and S. Stromme, {\it Bott's formula
and enumerative geometry}, J. Amer Math Soc {\bf 9} (1996) 175-193.

\bibitem{Gath} A. Gathmann, {\it Absolute and relative
Gromov-Witten invariants of very ample hypersurfaces},
math.AG/9908054.

\bibitem{Giv}
A. Givental, {\it Equivariant Gromov-Witten invariants},
Int Math Res Notices {\bf 13} (1996) 613--663.

\bibitem{GraP}
T. Graber and R. Pandharipande, {\it Localization of virtual classes}, Invent
Math {\bf 135} (1999), 487--518.

\bibitem{Kon}
M. Kontsevich, {\it Enumeration of rational curves with 
torus actions}, In: The moduli space of curves, (R. Dijkgraaf, C. Faber 
and G. van der Geer eds.) Prog in Math {\bf 129}, Birkh\"auser, Boston
(1995) 335-368.

\bibitem{KonM}
M. Kontsevich and Y. Manin,{\it Gromov-Witten classes, quantum
  cohomology, and enumerative geometry}, Comm Math Physics
{\bf 164} (1994), 525--562.

\bibitem{Lee}
Y.P. Lee, {\it Quantum Lefschetz hyperplane theorem},
math.AG/0003128.

\bibitem{LiT} J. Li and G. Tian, {\it Virtual moduli cycles and
Gromov-Witten invariants of algebraic varieties}, J Amer Math Soc
{\bf 11} (1998) 119-174.

\bibitem{LLY}
B. Lian, K. Liu, and S.T. Yau, {\it Mirror principle I}, Asian J
Math {\bf 1} (1997) 729--763. {\it Mirror principle II}, Asian
J Math {\bf 3} (1999) 109-146. {\it Mirror principle III}, Asian J of Math
{\bf 3} (1999) 771-800. {\it Mirror principle IV}, math.AG/0007104.

\bibitem{Ran} Z. Ran, {\it Enumerative geometry of singular plane
curves}, Invent Math {\bf 97} (1989) 447-465.

\end{thebibliography}
\end{document}